\documentclass[11pt]{article}
\usepackage{mathrsfs}
\usepackage{latexsym,lineno}
\usepackage{epsfig}
\usepackage{color}
\usepackage{amsmath}\usepackage{fleqn}\usepackage{verbatim}\usepackage{epsf}
\usepackage{amsthm}\usepackage{graphicx, float}\usepackage{graphicx}
\usepackage{amsfonts}\usepackage{amssymb}\usepackage{graphpap}
\usepackage{epic}\usepackage{curves}

\newcommand{\be}{\begin{equation}}
\newcommand{\ee}{\end{equation}}
\newcommand{\benum}{\begin{enumerate}}
\newcommand{\eenum}{\end{enumerate}}
\newcommand{\bit}{\begin{itemize}}
\newcommand{\eit}{\end{itemize}}

\newtheorem{prop}{Proposition}[section]

\usepackage{graphicx}

\begin{document}
\def\s{\subseteq}
\def\n{\noindent}
\def\se{\setminus}
\def\dia{\diamondsuit}
\def\la{\langle}
\def\ra{\rangle}


\title{On extremal multiplicative Zagreb indices of trees with given number of vertices of maximum degree}

\author{Shaohui Wang$^a$, ~Chunxiang Wang$^{b,}$\footnote{  Authors'email address: S. Wang (e-mail:
shaohuiwang@yahoo.com; swang@adelphi.edu), C. Wang(e-mail: wcxiang@mailccnu.edu.cn), L. Chen(e-mail: 409137414@qq.com)
J.B. Liu (e-mail:  liujiabaoad@163.com).}~, ~Lin Chen$^b$, Jia-Bao Liu$^c$   \\
\small\emph {a. Department of Mathematics and Computer Science, Adelphi University, Garden City, NY 11550, USA}\\
\small\emph {b. School of Mathematics and Statistics, Central China Normal University, Wuhan
430079, P.R. China} \\
\small\emph {c. School of Mathematics and Physics, Anhui Jianzhu
University, Hefei 230601, P.R. China} }
\date{}
\maketitle

\begin{abstract}

The first multiplicative Zagreb index of a graph $G$ is the product of the square of every vertex degree, while the second multiplicative Zagreb index is the product of the products of degrees of pairs of adjacent vertices.   In this paper, we explore the trees in terms of given  number of vertices of maximum degree. The maximum and minimum values of $\prod_1(G) $ and $\prod_2(G) $ of trees with arbitrary number of maximum degree are provided.  In addition, the corresponding extremal graphs are characterized.\\
Accepted by Discrete Applied Mathematics.

\end{abstract}

\section{Introduction}
Throughout this paper, we consider simple, connected and undirected graphs.
Denote a graph by $G = (V, E)$, where $V = V (G)$ is called vertex
set and $E = E(G)$ is called edge set. For a vertex $v \in V(G)$, the neighborhood of  $v $ is the set $N(v) = N_G(v) = \{w \in V(G), vw \in E(G)\}$, and $d_G(v)$ (or   $d(v)$)  denotes the degree of $v$ with $ d_{G}(v) = |N(v)|$. $n_i$ is the number of vertices of degree $i \geq 0$.
 If a graph $G$  contains $n$
vertices and $n-1$ edges, then $G$ is called a tree. For a
vertex $v \in V (T )$ with $2 \leq d_T(v) \leq \Delta(T)-1$,  its edge rotating capacity is defined to be $d_T (v)-1$. The total edge rotating capacity of a tree $T$ is equal to the sum of the
edge rotating capacities of its vertices that satisfy the condition $2 \leq d_T(v) \leq \Delta(T)-1$.  As usual,    denote $P_n$  by the path   on $n$ vertices. The maximum vertex degree in the graph $G$  is  denoted by $\Delta(G)$.

The degree sequence of G is a sequence of positive integers $\pi=(d_1, d_2, \cdots, d_n)$  if $d_i=d_G(v) (i = 1,\cdots, n)$ holds for $v \in V(G)$. In this work, we assign an order of the vertex degrees as non-increasing, i.e., $d_1\geq d_2 \geq \cdots \geq d_n$. In addition, a sequence $\pi =(d_1, d_2, \cdots, d_n)$ is called a tree degree sequence if there exists a tree T such that $\pi$ is its degree sequence.  Furthermore, it is well known that the sequence $\pi =(d_1, d_2, \cdots, d_n)$  is a degree sequence of a  tree with $n$ vertices if and only if
$$\sum_{i=1}^{n}d_i=2(n-1).$$

In the interdisplinary of mathemactics, chemistry and physics,
 molecular invariants/descriptors could be useful for the study of quantitative structure-property relationships (QSPR) and quantitative structure-activity relationships (QSAR)
 and
for the descriptive presentations of biological and chemical
properties, such as boiling and melting points, toxicity,
physico-chemical, and biological
properties~\cite{Gutman1996,Liu2015,LiuP2015,LiuPX2015,x001,0001,0002}. One class of the oldest topological molecular descriptors
are named as Zagreb indices~\cite{Gutman1972}, which are literal
quantities in an expected formulas for the total
$\pi$-electron energy of conjugated molecules as follows.
\begin{eqnarray} \nonumber
M_1(G) = \sum_{u \in V(G)} d(u)^2
~\text{ and}
~ M_2(G) = \sum_{uv \in E(G)} d(u)d(v).
\end{eqnarray}

Based on the successful considerations on
these applications of Zagreb indices~\cite{Gutman2014},
Todeschini et al.(2010)~\cite{RT20101,RT20102,Wang2015} presented the
following multiplicative variants of molecular structure
descriptors:
\begin{eqnarray} \nonumber
 \prod_1(G) = \prod_{u \in V(G)} d(u)^2 ~
\text{ and}\;\;
 \prod_2(G) = \prod_{uv \in E(G)} d(u)d(v) = \prod_{u \in V(G)}
d(u)^{d(u)}.
\end{eqnarray}

  Recently, there are lots of articles
 explored multiplicative Zagreb indices in the interdisplinary of chemistry and
mathematics~\cite{Hu2005,Li2008,shi2015,BF2014,SM2014,Xu2014,WangJ2015}.
Iranmanesh et al.~\cite{Iranmanesh20102}  explored first and second multiplicative Zagreb indices for a class of chemical dendrimers.
 Xu and Hua~\cite{Xu20102} provided an unified approach to
characterize extremal maximal and minimal trees, unicyclic
graphs and bicyclic graphs regarding to multiplicative Zagreb
indices, respectively.
Wang and Wei~\cite{Wang2015} gave the maximum and minimum indices of these indices in $k$-trees, and the corresponding extreme graphs are provided.
Liu and Zhang [14] investigated some sharp upper bounds for
$\prod_1$-index and $\prod_2$-index in terms of graph parameters
such as an order, a size and a radius~\cite{Liuz20102}.
Kazemi~\cite{Ramin2016} studied the bounds for
the moments and the probability generating function of these
indices in a randomly chosen molecular graph with tree structure
of order $n$.
Borovi\'canin et al.~\cite{Borov2016} introduced upper
bounds on Zagreb indices of trees, and a lower bound for the
first Zagreb index of trees with a given domination number is
determined and the extremal trees are characterized as well.
 Borovi$\acute{c}$anin and Lampert\cite{Bojana2015} provided the maximum and minimum
Zagreb indices of trees with given
number of vertices of maximum degree.

Motivated by  above results, in this paper we further
investigate the multiplicative Zagreb indices of trees with arbitrary  number of vertices of maximum degree.  The maximum and minimum values of $\prod_1(G) $ and $\prod_2(G) $ of trees with arbitrary number of maximum degree are provided.  In addition, the
corresponding extreme graphs are charaterized.
Our results extends and enriches  some known conclusions
obtained by \cite{Bojana2015}.

\section{Preliminaries}

It is known that each tree has at least two minimum degree vertices, named as pendent vertices, and
some maximum degree vertices.  It is natural to consider the trees with arbitrary  number of
maximum degree vertices.

Let $\mathcal{T}_{n, k}$ be the class of trees with $n$ vertices, in which there exist $k$ vertices having the maximum degree with $n > k > 0$.  Note that the path $P_n$ is the unique element of $\mathcal{T}_{n,n-2}$. So, in the following we consider the class $\mathcal{T}_{n,k}$ with $k\leq n-3$.

We first introduce several facts and tools, which are important in the proofs of following sections.

\begin{prop}\cite{Bojana2015} If $T\in \mathcal{T}_{n,k}$ is a tree with $k$  vertices of maximum  degree $\Delta$, then $\Delta \leq \lfloor   \frac{n-2}{k}\rfloor +1$.
\end{prop}

By the routine calculations, one can derive the following propositions.
\begin{prop}
Let $f(x) = \frac{x}{x+m}$ be a function with $m > 0$. Then $f(x)$ is increasing in $\mathbb{R}$.
\end{prop}

\begin{prop}
Let $g(x) = \frac{x^x}{(x+m)^{x+m}}$ be a function with $m > 0$. Then $g(x)$ is decreasing  in $\mathbb{R}$.
\end{prop}

Based on   the above algebraic tools, we are ready to provide the sharp upper and lower bounds of first multiplicative Zagreb index of such trees in section 3, and the sharp upper and lower bounds of second multiplicative Zagreb index of these trees in section 4. Some of notations and figures are used close to \cite{Bojana2015}.

\section{  The sharp upper and lower bounds of first mutiplicative Zagreb index on the trees
 }

In this section, we  obtain the  bounds of the first multiplicative Zagreb index in the class $\mathcal{T}_{n,k}$.

\subsection{The sharp upper bounds of $\prod_1$   on trees with given number of vertices of maximum degree}
The first multiplicative Zagreb index   of $\mathcal{T}_{n, k}$ can be routinely calculated   if the degree sequence is given.

\noindent {\bf  Lemma 3.1. } Let ${T}_{min}^1$ be a tree with minimal value of  first multiplicative Zagreb index in $\mathcal{T}_{n,k}$. Then $\Delta ({T}_{min}^1)  =\lfloor \frac{n-2}{k}\rfloor +1$.

\begin{proof} Let $\Delta $ be the maximum vertex degree in the tree $T_{min}^1$. By  Proposition 2.1, we have $\Delta\leq \lfloor   \frac{n-2}{k}\rfloor +1$.  Denote $ \Delta_{max}=\lfloor   \frac{n-2}{k}\rfloor +1$ and $n-2 = k \lfloor   \frac{n-2}{k}\rfloor +r$, where $0\leq r <k$.

Firstly, we assume that  $\Delta<\Delta_{max}$.
Denote $V({T}_{min}^1)$ =$\{v_1, \cdots, v_n\}$ by    the degree sequence: $\pi=(d_1, d_2, \cdots , d_n)$. Then
$$\Delta=d_1=\cdots=d_k=\Delta_{max}-t, \; t>0.$$
Note that  \begin{eqnarray}
\sum_{i=1}^{\Delta}n_i=n
\end{eqnarray}
 and
 \begin{eqnarray}
\sum_{i=1}^{\Delta} i n_i=2(n-1).
\end{eqnarray}
By the relations (1) and (2),  we can obtain that
\begin{eqnarray}
n_1\geq 2+(\Delta -2)n_{\Delta}=2+k(\Delta -2).
\end{eqnarray}
Let $n_1=2+k(\Delta -2)+n_1'$ with $n_1'\geq 0$. Using (1) we obtain
 $$n=2+k(\Delta -2)+n_1'+\sum_{i=2}^{\Delta-1} n_i+k,$$
 which implies   that
\begin{eqnarray}
 \sum_{i=2}^{\Delta-1} n_i +n_1'=r+kt.
\end{eqnarray}
Also,  by the relation (2),  one can calculate that
 \begin{eqnarray}
 \sum_{i=2}^{\Delta-1} in_i +n_1'=2(r+kt).
\end{eqnarray}
By subtracting the relations (4) and (5),  we have
 \begin{eqnarray}
 \sum_{i=2}^{\Delta-1} (i-1)n_i =r+kt\geq kt \geq k.
\end{eqnarray}

Then the total edge rotating capacity of this tree is greater than or equal to $k$.

  Let $v_i\in V({T}_{min} ^1)$ ($k<i\leq n)$ be a vertex that has positive edge rotating capacity and let
  $d_i$  be its degree with $2\leq d_i \leq \Delta -1$.   Define a tree ${T}_1$ with the vertex degree sequence
  $\pi_1=(d_1^1, d_2^1, \cdots, d_n^1)$ such that $d_1^1=\Delta +1$, $d_i^1=d_i -1$ and $d_j^1=d_j(j\in \{2,\cdots, n\}$, $j\ne i )$.
By the definition of $\prod_1$, we have
 \begin{eqnarray}
 \nonumber  \frac{ \prod_1({T}_1)}{\prod_1({T}_{min}^1)} &&=\frac{(d_1+1)^2d_2^2\cdots(d_i-1)^2\cdots d_n^2}{d_1^2d_2^2 \dots d_i^2 \dots d_n^2}
\\  \nonumber &&= (\frac{d_1+1}{d_1})^2(\frac{d_i-1}{d_i})^2
\\  \nonumber &&= (\frac{\frac{(d_i-1)}{(d_i-1)+1} }{\frac{d_1}{d_1+1}})^2~ \;\;\;\;(\text{by Proposition} \;2.2)
  \\
&&<1.
\end{eqnarray}

 Therefore, $\prod_1({T}_1) <\prod_1({T}_{min}^1)$. Note that ${T}_1$ has the maximum vertex degree equal to $\Delta+1$, then $T_1\notin \mathcal{T}_{n,k}$. Since ${T}_{min}^1$ has $k-1$ vertices of degree $\Delta$, and the total positive edge rotating capacity of $d_2$, $d_3$, \dots, $d_{\Delta-1}$ is at least $k-1$, by the relation $(6)$. Then we can conclude the following statements.

 We recursively proceed the above described transformations of the tree ${T}_{min}^1$ for $k-1$ times on each vertex of maximum degree $\Delta$. In every step,  we could define a tree ${T}^l$ with  degree sequence $\pi_l=(d_1^l, \dots ,d_n^l)$ such that $d_l^l=\Delta+1$, $d_i^l=d_i^{l-1}-1$ and $d_j^l=d_j^l-1(j\in \{1, \cdots, n\}$, $j\ne i, \;$l$ )$, where $l=2,\cdots,k$ and $d_i^{l-1}$ is the degree of any vertex $v_i\in V(T^{l-1})$ with $k<i\leq n$, which has positive edge rotating capacity (this vertex exists, because the total edge rotating capacity of the tree $T^{l-1}$ is at least $k-l+1)$. It is natural to see that after some described transformations,  we obtain a tree whose degrees $d_{k+1} , \cdots , d_n$ are in an  increasing order. Except for these steps, every such transformation strictly decreases the first multiplicative Zagreb index. Finally, we could get trees  ${T}^k, {T}^{k-1}, \cdots, {T}^2  \in \mathcal{T}_{n,k}$, in which they have the maximum vertex degree equal to $\Delta+1=(\Delta_{max}-t)+1$ and satisfy the conditions $\prod_1({T}^k)<\prod_1({T}^{k-1})< \cdots <\prod_1({T}^1)<\prod_1({T}_{min}^1)$. Considering that ${T}^k$ has $k$ vertices of degree $\Delta+1$, all of them are vertices of maximal degree.

 Since $\prod_1({T}^k) <\prod_1({T}_{min}^1)$, it contradicts the fact that ${T}_{min}^1$ has the minimal first multiplicative Zagreb index in $\mathcal{T}_{n,k}$. Thus, $t=0$ and it could be considered as $\Delta=\Delta_{max}=\lfloor \frac{n-2}{k}\rfloor +1$.

 This completes the proof of  $\Delta=\Delta_{max}=\lfloor \frac{n-2}{k}\rfloor +1$.
\end{proof}

Based on the above proof, the following remark is immediate.

 \noindent {\bf  Remark 3.1. } The statement of Lemma 3.1   holds for $k=n-2$, i.e., $ T_{n, n-2}  =P_n$.

The next theorem provides the sharp lower bound of $\prod_1(G)$ and charaterizes the exremal graphs achieving such lower bound.

\noindent {\bf  Theorem 3.1. } Let ${T}\in \mathcal{T}_{n,k}$, where $1\leq k\leq \frac{n}{2}-1$. Then
$$\prod_1({T})\geq\Delta^{2k}(\Delta-1)^{2p}\mu^{2},$$
where the equality holds if and only if its   degree sequence  is $(\underbrace{\Delta,\cdots,\Delta}_{k}, \underbrace{\Delta-1,\cdots,\Delta-1}_p, \mu, \underbrace{1,1,\cdots,1}_{n-k-p-1})$ with $\Delta=\lfloor \frac{n-2}{k}\rfloor +1$, $p=\lfloor \frac{n-2-k(\Delta-1)}{\Delta-2}\rfloor $ and $\mu=n-1-k(\Delta-1)-p(\Delta-2)$.

\begin{proof}
 Let $\pi=(d_1,d_2,\cdots,d_n)$ be the vertex degree sequence of a tree $T_{min}^1$ with minimal first multiplicative Zagreb index in $T_{n,k}$. By Lemma 3.1 we have that $d_1=d_2=\cdots=d_k=\Delta={\lfloor \frac{n-2}{k}\rfloor+1}$. Since $\Delta-1=\lfloor \frac{n-2}{k}\rfloor$, then one can obtain that $\Delta-1$ is the integer section of $\frac{n-2}{k}$. Based on the previous lemma, let $n-2=k(\Delta-1)+r$, where $0\leq r<k$. From the relation (3), it follows that $n_1\geq k(\Delta-2)+2=n-k-r$, and $n_1$ has at least $n-k-r$ vertices with one degree. Therefore, $d_n=d_{n-1}=\dots=d_{k+r+1}=1$.

Note that $n_{\Delta}=k$ and $n_1=k(\Delta+2)+2+n_1'$ with $n_1'\geq 0$. Combining with the relations (1) and (2), and $t=0$ in the relations (5) and (6), we obtain that

\begin{eqnarray}
 \sum_{i=2}^{\Delta-1}n_i+n_1'=r
\end{eqnarray}
and
\begin{eqnarray}
\sum_{i=2}^{\Delta-1}(i-1)n_i=r.
\end{eqnarray}

Since $n_i\geq 0$ with $i=2,3,\cdots,\Delta-1$, from the relation (9) it follows that $p=n_{\Delta-1}\leq\frac{r}{\Delta-2}$. Because $p$ is a non-negative integer number,  we obtain $p\leq \lfloor \frac{r}{\Delta-2}\rfloor$.

Now,  suppose that $p< \lfloor \frac{r}{\Delta-2}\rfloor$ and  $\lfloor \frac{r}{\Delta-2}\rfloor=t'$.  Let
 \begin{eqnarray}
 r=n_2 + 2n_3 +\cdots + (\Delta-3)n_{\Delta-2}+(\Delta-2)p=t'(\Delta-2) + y,\end{eqnarray}
where $0\leq y<(\Delta-2)$.  Then $\sum_{i=2}^{\Delta-2}(i-1)n_i\geq {\Delta-2}$.

So there exist $n_i$ and $n_j(2\leq i< j\leq {\Delta-2})$, where $n_i\ne 0$ and $n_j\ne 0$ or $n_i\geq 2$~(where $2\leq i\leq \Delta-2)$ and the equality (9) is satisfied. Furthermore, since $\pi=(\underbrace{\Delta,\Delta,\cdots,\Delta}_k,d_{k+1},\cdots,d_{k+r},1,1,\cdots,1)$, then there exist numbers $d_{k+j_1}$ and $d_{k+i_1} (1\leq j_1< i_1\leq r)$  such that $d_{k+j_1}=j> d_{k+i_1}=i ~(\text{or}~ d_{k+j_1}=d_{k+i_1}=i, ~\text{if}~ n_i\geq 2)$.

Let $\pi'=(d_1',d_2',\cdots,d_n')$ be a sequence of positive integers such that   $d_{k+j_1}'=d_{k+j_1}+1=j+1 $ and $d_{k+i_1}'=d_{k+i_1}+1=i-1 $, where $d_u'=d_u$ (for $u\ne {k+j_1}, u\ne {k+i_1})$. Clearly, $\sum_{i=1}^{n}d_i'=2n-2$, and $\pi'$ is the vertex degree sequence of a tree $T'$. Also,
 \begin{eqnarray}
 \nonumber \frac{ \prod_1({T}')}{\prod_1({T}_{min}^1)} &&=\frac{(j+1)^2(i-1)^2}{j^2 i^2}
\\  \nonumber &&= (\frac{\frac{i-1}{(i-1)+1} }{\frac{j}{j+1}})^2 ~\;\;\;\;(\text{by Proposition 2.3})
 \\
&&<1.
\end{eqnarray}

Therefore, $\prod_1({T}')<\prod_1({T}_{min}^1)$. This contradicts the choice of ${T}_{min}^1$ which has the minimal first multiplicative Zagreb index in the class $\mathcal{T}_{n,k}$.  Hence we  conclude that $p=n_{\Delta-1}=\lfloor \frac{r}{\Delta-2}\rfloor$.

Next, by the relation (10) we obtain that
$y$=$n_2$+$2n_3$+ $\cdots$ +$(\Delta-3)n_{\Delta-2}$~ (for $0\leq y<
{\Delta-2}$),  and
\begin{eqnarray}
 y=r-p(\Delta-2)\leq {\Delta-3}.
\end{eqnarray}

According to the analysis of the relation (10), it can  be proved that   $n_\mu=1$ (where $2\leq \mu\leq {\Delta-2}$). Then $\mu= y+1={r-p(\Delta-2)+1}$, i.e., $\mu=n-1-k(\Delta-1)-p(\Delta-2)$, and $n_i=0$, for $i\ne \mu$, where $2\leq i\leq{\Delta-2}$, since in the opposite case we can  construct a tree whose $\prod_1$ is smaller than $\prod_1({T}_{min}^1)$ again.

  Hence, the tree ${T}_{min}^1$ has minimum first multiplicative Zagreb index.   So we conclude that $n_\mu=1$, and $\mu=r-p(\Delta-2)+1 $, namely, $\mu={n-1-k(\Delta-1)-p(\Delta-2)}$.

Therefore, the tree ${T}_{min}^1$ with minimum first multiplicative Zagreb index in the class $\mathcal{T}_{n,k}$ has the vertex degree sequence $\pi=(\underbrace{\Delta,\cdots,\Delta}_{k}, \underbrace{\Delta-1,\cdots,\Delta-1}_p, \mu, \underbrace{1,1,\cdots,1}_{n-k-p-1})$.
\end{proof}

The first multiplicative Zagreb index of the tree $\mathcal{T}_{min}^1$ can now be routinely calculated.

 \noindent{\bf  Remark 3.2. } The statement of Theorem 3.1 also holds for $k=n-2$, i.e., $ T_{n,k}=P_n$.

\subsection{The sharp  lower bound of $\prod_1$  on trees with given number of vertices of maximum degree
 }

     In the following theorem we will describe the trees that have the maximal first multiplicative Zagreb index in the class $\mathcal{T}_{n,k}$.

 \noindent {\bf  Lemma 3.2. } Let ${T}_{max}^1$ be a tree with maximal first multiplicative Zagreb index in the class $\mathcal{T}_{n,k}$, where $1\leq k\leq {\frac{n}{2}-1}$. Then its maximum   degree $\Delta$ equals to 3.

 \begin{proof}

We first assume that $\Delta \geq 4$ and $u$ is a vertex of maximum degree $\Delta$ in ${T}_{max}^1$.

\begin{figure}[htbp]
    \centering
    \includegraphics[width=4in]{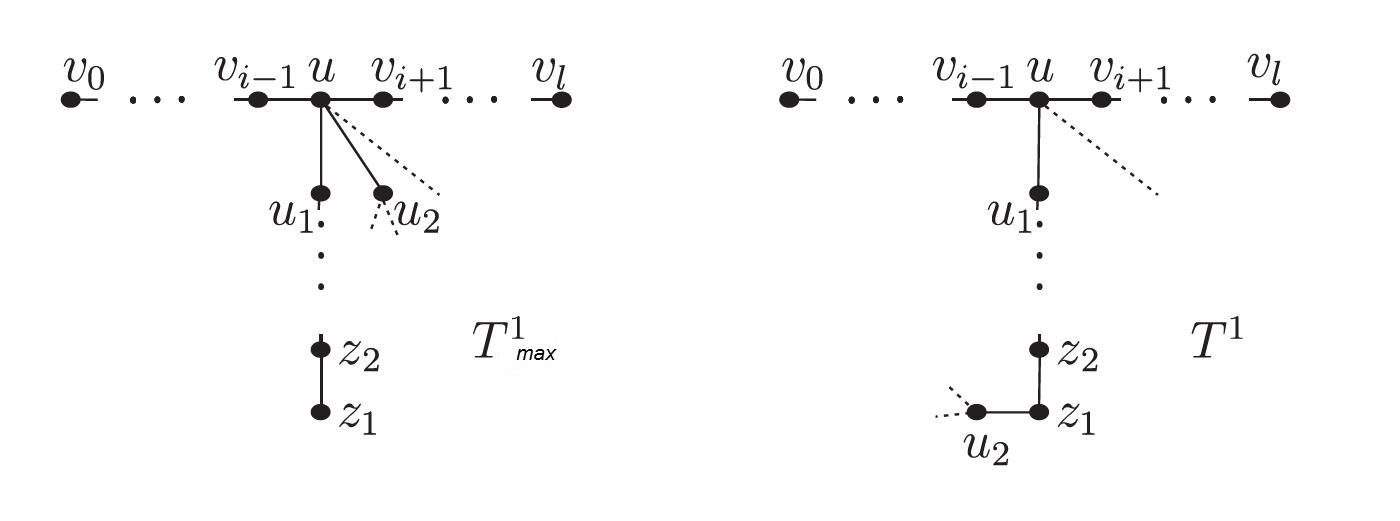}
    \caption{ The graphs ${T}_{max}^1$ and ${T}^1$ in Lemma 3.2.}
    \label{fig: te}
\end{figure}

 Denote $p=v_0v_1...v_{i-1}u(=v_i)v_{i+1}...v_l$ by the longest path in ${T}_{max}^1$ that contains $u$.
 Also, let $v_{i-1}$, $v_{i+1}$, $u_1$, $u_2$, ...$u_{\Delta-2}$ be the vertices adjacent to $u$ in ${T}_{max}^1$, and $z_1$ be a pendent vertex connected to $u$ via $u_1$ (it is possible that $z_1\equiv u_1$)(Figure 1). Now we define a tree ${T}^1$ such that
\begin{eqnarray}
  {T}^1={T}_{max}^1 - uu_2 + u_2 z_1.
\end{eqnarray}
Then
\begin{eqnarray}
  \frac{\prod_1({T}^1)}{\prod_1({T}_{max}^1)}  = \frac{(\Delta-1)^2 2^2}{\Delta^2 1^2}
=  ({\frac{\Delta-1}{\Delta}})^2 4
>1,
\end{eqnarray}
that is, $\prod_1({T}^1)>\prod_1({T}_{max}^1)$. Obviously, the tree ${T}^1$ has $k-1$ vertices of degree $\Delta$.

 Analogously, we  apply similar transformations described in the relation (13) on every vertex $u$ of degree $\Delta$. In each step from a tree $T^i$, we   obtain a tree ${T}^{i+1}$ with $1\leq i\leq k-1$ that has greater first multiplicative Zagreb index than its predecessor. After $k$ repetitions of these transformations, we arrive at the tree $T^k$, which has $k$ vertices having the maximum degree $\Delta-1$. Clearly, ${T}^k\in \mathcal{T}_{n,k}$ and $\prod_1({T}^k)>\prod_1({T}_{max}^1)$.

 Considering the tree  maximizing $\prod_1$ in the class $\mathcal{T}_{n,k}$, it contradicts the choice of ${T}_{max}^1$. Then ${T}_{max}^1$ has the maximum degree of $3$, this completes our proof.
\end{proof}

\noindent {\bf  Theorem 3.2. } Let ${T}\in \mathcal{T}_{n,k}$ with $1\leq k\leq \frac{n}{2}-1$. Then
  $$\prod_1({T})\leq 9^{k}4^{(n-2k-2)}=(\frac{9}{16})^{k}4^{n-2},$$
 where the equality holds if and only if  ${T}$ has   degree sequence
 $\pi=(\underbrace{3,3,\cdots,3}_k, \underbrace{2,2,\cdots,2}_{n-2k-2}, \underbrace{1,1,\cdots,1}_{k+2})$.

\begin{proof}
 Let $T_{max}^1$ be a tree with maximum $\prod_1$ in the class $\mathcal{T}_{n,k}$. By Lemma 3.2, we obtain $\Delta=3$. So the vertex degree sequence of this tree is $\pi=(\underbrace{3,3,\cdots,3}_k, \underbrace{2,2\cdots,2}_{n_2}, \underbrace{1,1,\cdots,1}_{n_1})$ with $k\leq \frac{n}{2}-1$. Hence, applying to the equality (1) we obtain
\begin{eqnarray}
  n_1   + 2n_2 + 3k = 2(n-1) = 2(n_1+n_2+k)-2.
\end{eqnarray}
Furthermore, since $n_3=\Delta=k$, then
\begin{eqnarray}
  n_1 + 2n_2 + 3k = 2(n-1) = 2n_1 + 2n_2 + 2k -2.
\end{eqnarray}
 According to the relation (16), we can conclude that $n_1$=$k$+2 and $n_2$=$n$-$(n_1+n_3)$=$n$-$2k$-2.
Hence, $$\pi=(\underbrace{3,3,\cdots,3}_k, \underbrace{2,2,\dots,2}_{n-2k-2}, \underbrace{1,1,\cdots,1}_{k+2}),$$
 and $$\prod_1({T}_{max}^1)=d_1^2d_2^2\cdots d_n^2={\underbrace{3^2 3^2\cdots 3^2}_k}{\underbrace{2^22^2\cdots 2^2}_{n-2k-2}}{\underbrace{1^21^2 \cdots 1^2}_{k+2}} =9^{k}4^{n-2k-2} =(\frac{9}{16})^k4^{n-2}.$$
Therefore, this completes our proof.
\end{proof}

\noindent {\bf  Remark 3.3. } Let ${T}^k$ be a tree with maximal first multiplicative Zagreb index in the class $\mathcal{T}_{n,k}$, and $T^p\in \mathcal{T}_{n,p}$, where $1\leq {k, p}\leq {\frac{n}{2}-1}$. If $k> p$, then $\prod_1({T}^k)< \prod_1({T}^p)$.

\begin{proof} We can directly obtain by Theorem 3.2.
\end{proof}

\section{ The sharp upper and lower bounds of second mutiplicative Zagreb index on the trees  }

In this section, we first characterize the trees with maximal second multiplicative Zagreb index in the class $\mathcal{T}_{n,k}$, and let $\pi=(d_1,d_2,\cdots,d_n)$ be a degree sequence of a tree. To prove our main result in this section,  we provide  a lemma in each subsection and use it to deduce our theorems.

\subsection{The sharp  lower bound of $\prod_2$ on trees with given number of vertices of maximum degree}
\noindent{\bf  Lemma 4.1. } Let $ T_{min}^2$ be a tree with minimal second multiplicative Zagreb index in the class $\mathcal{T}_{n,k}$. Then, $\Delta(T_{min}^2)=3$.

\begin{proof}
Suppose that $\Delta \geq 4$ and $u$ is a vertex of maximum degree $\Delta$ in $\mathcal{T}_{min}^2$.

\begin{figure}[htbp]
    \centering
    \includegraphics[width=4in]{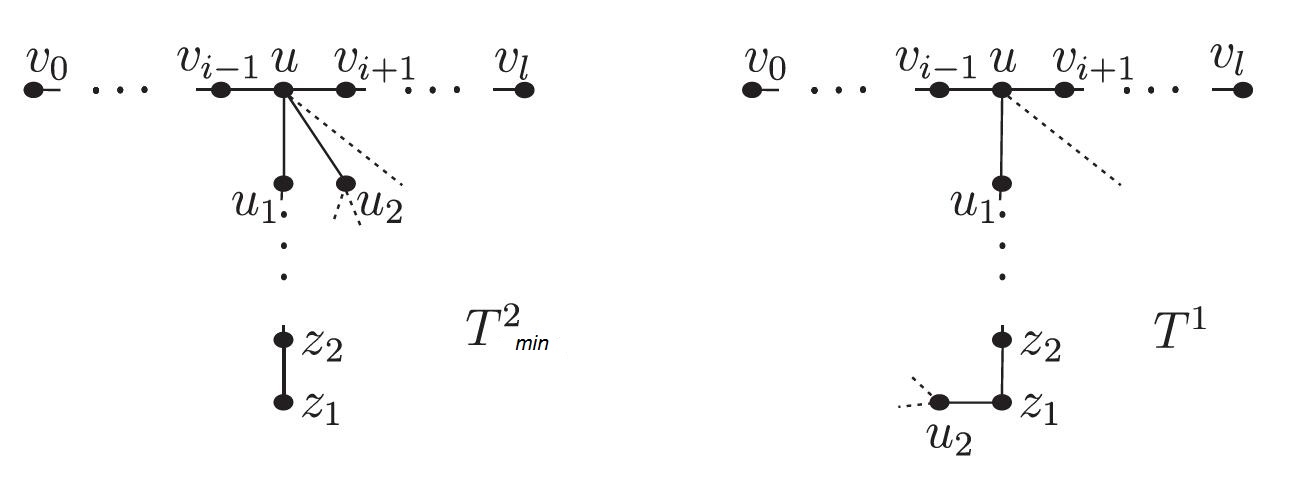}
    \caption{ The graphs ${T}_{min}^2$ and ${T}^1$ in Lemma 4.1.}
    \label{fig: te}
\end{figure}

Let $p=v_0v_1...v_{i-1}u(=v_i)v_{i+1}\cdots v_l$ be the longest path in ${T}_{min}^2$ that contains $u$. Also, let $v_{i-1}$, $v_{i+1}$, $u_1$, $u_2$, $\cdots$,$u_{\Delta-2}$ be the vertices adjacent to $u$ in ${T}_{min}^2$, and $z_1$ a pendent vertex connected to $u$ via $u_1$ (it is possible that $z_1\equiv u_1$)(Figure 2). Likewise, we define a tree ${T}^1$ such that
\begin{eqnarray}
 {T}^1= {T}_{min}^2 - uu_2 + u_2z_1.
\end{eqnarray}
Then
\begin{eqnarray}
 \nonumber   \frac{\prod_2({T}^1)}{\prod_2({T}_{min}^2)} &&= \frac{(d_u-1)^{d_u-1} d_{z_1}^{d_{z_1}}}{d_u^{d_u} d_{z_1}^{d_{z_1}}}
\\  \nonumber &&=   \frac{(d_u-1)^{d_u-1}{2^2}}{d_u^{d_u}{1^1}}
\\  \nonumber &&=4 \frac{(\Delta-1)^{\Delta-1}}{\Delta^\Delta} ~\;\;\;\;(\text{by Proposition 2.3 and}~
  \Delta-1  \geq 3)
   \\ &&< 4\frac{3^3}{(3+1)^{3+1}}< 1.
 \end{eqnarray}

Therefore, $\prod_2({T}^1)< \prod_2({T}_{min}^2)$.

In the same way, we   use these transformations described in the relation (16) on every vertex $v$ of degree $\Delta$. Now, we   transform ${T}_{min}^2$ with the remaining  $k-1$ vertices of the degree $\Delta$ into $\Delta-1$. After $k$ repetitions of the transformation we arrive at the trees ${T}^2$, ${T}^3$,$\cdots$,${T}^k$, in which $k$ vertices have maximum degree $\Delta-1$. Obviously, ${T}^k\in \mathcal{T}_{n,k}$ and
 $\prod_2({T}^k)<
                     \prod_2({T}^{k-1})<\cdots<\prod_2({T}_{min}^2)$.

Thus, ${T}^k$ has $k$   maximum degree vertices. This contradicts the choice of ${T}_{min}^2$ as the tree that minimizes $\prod_2$ in the class $\mathcal{T}_{n,k}$.
Therefore, $\Delta({T}_{min}^2)=3$.
\end{proof}

\noindent{\bf  Theorem 4.1. } Let ${T}\in \mathcal{T}_{n,k}$, where $1\leq k\leq \frac{n}{2}-1$. Then
$$\prod_2({T})\geq {(3^3)^k}{(2^2)}^{n-2k-2}=({\frac{27}{16}})^k 4^{n-2},$$
where the equality holds if and only if  ${T}$ has degree sequence
 $\pi=(\underbrace{3,3,\cdots,3}_k, \underbrace{2,2,\cdots,2}_{n-2k-2}, \underbrace{1,1,\cdots,1}_{k+2})$.

\begin{proof}
 Let ${T}_{min}^2$ be a tree with minimal $\prod_2(G)$ in the class $\mathcal{T}_{n,k}$. According to Lemma 4.1, we obtain $\Delta({T}_{min}^2)=3$.
Thus,
$$\pi=(\underbrace{3,3,\cdots,3}_k, \underbrace{2,2\cdots,2}_{n_2}, \underbrace{1,1,\cdots,1}_{n_1}).$$
Clearly, $n_1+2n_2+3k=2(n-1)=2n_1+2n_2+2k-2.$ Then $n_1=k+2$, $n_2=n-n_1-n_3=n-2k-2$.
Therefore, $\prod_2({T}_{min}^2)={(3^3)^k}{(2^2)^{n_2}}{(1^1)^{n_1}}={27}^k4^{n-2k-2}=(\frac{27}{16})^k 4^{n-2}$.
\end{proof}

\noindent {\bf  Remark 4.1. } Let ${T}^k$ be a tree with minimal second multiplicative Zagreb index in $\mathcal{T}_{n,k}$, and ${T}^p$ a tree with minimal second multiplicative Zagreb index in the class $\mathcal{T}_{n,p}$, where $1\leq {k, p}\leq {\frac{n}{2}-1}$. If $k> p$, then $\prod_2({T}^k)> \prod_2({T}^p)$.

\begin{proof} According to Theorem 4.1 we can deduce that the proposition is correct.
\end{proof}

\subsection{The sharp  upper bound of $\prod_2$ on   trees with given number of vertices of maximum degree
 }
\noindent {\bf  Lemma 4.2. } Let ${T}_{max}^2$ be a tree with maximum second multiplicative Zagreb index in the class $\mathcal{T}_{n,k}$, then its maximum vertex degree $\Delta({T}_{max}^2)=\lfloor \frac{n-2}{k}\rfloor +1$.

\begin{proof} Let $\Delta=\Delta({T}_{max}^2)$  be the maximum vertex degree of a tree ${T}_{max}^2$. By Proposition 2.1, $\Delta \leq \lfloor \frac{n-2}{k} \rfloor +1$. By the similar proof of Lemma 3.1, we can conclude that the tree ${T}^1$.  Then,
\begin{eqnarray}
  \nonumber   \frac{\prod_2({T}^1)}{\prod_2({T}_{max}^2)}  &&=  \frac{(\Delta+1)^{\Delta+1} {(d_i-1)^{d_i-1}}}{\Delta^{\Delta} d_i^{d_i}}
\\  \nonumber &&=\frac{\frac{(d_i-1)^{(d_i-1)}}{((d_i-1)+1)^{((d_i-1)+1)}}} {\frac{\Delta^{\Delta}} {(\Delta+1)^{\Delta+1}}}~\;\;\;\;  (\text{by Proposition  2.3 and} ~ 2\leq d_i\leq \Delta-1<\Delta  )
\\    &&>1.
 \end{eqnarray}

Therefore, $\prod_2({T}^1)> \prod_2({T}_{max}^2)$. We can repeat   remaining described transformations of the tree ${T}_{max}^2$ on every of degree $\Delta$. Therefore, we can conclude that this contradicts the fact that ${T}_{max}^2$ has the maximum second multiplicative Zagreb index in the class $\mathcal{T}_{n,k}$.

This proves that $\Delta({T}_{max}^2)=\lfloor \frac{n-2}{k}\rfloor +1$.
\end{proof}

\noindent {\bf  Theorem 4.2. } Let ${T}\in \mathcal{T}_{n,k}$, where $1\leq k\leq
\frac{n}{2}-1$.  Then
 $$\prod_2( T )\leq \Delta^{k\Delta}(\Delta-1)^{p(\Delta-1)} \mu^\mu,$$
 where the equality holds if and only if the tree ${T}$ has a vertex degree sequence
$\pi=(\underbrace{\Delta,\Delta,\cdots,\Delta}_k, \\  \underbrace{(\Delta-1),(\Delta-1),\cdots,(\Delta-1)}_p, \mu,  \underbrace{1,1,\cdots,1}_{n-k-p-1})$,
for $\mu=n-1-k(\Delta-1)-p(\Delta-2)$ and $p=\lfloor \frac{n-2-k(\Delta-1)}{\Delta-2}\rfloor.$

\begin{proof}
Let ${T}_{max}^2$ be a tree with maximum $\prod_2$ in the class $\mathcal{T}_{n,k}$. By Lemma 4.2, we have that $\Delta({T}_{max}^2)=\lfloor \frac{n-2}{k}\rfloor +1$.
Let
$p =n_{\Delta-2}$, then $p=\lfloor \frac{r}{\Delta-2}\rfloor. $
Otherwise, we obtain the tree ${T}'$ (using the method is analogous to the Theorem 3.1  which constructs the tree $T'$ )  such that $\pi'$ is the vertex degree sequence of a tree $T'$. Also, $\pi'=(d_1',d_2',\cdots,d_n')$ is a sequence of positive integers such that   $d_{k+j_1}'=d_{k+j_1}+1=j+1 $ and $d_{k+i_1}'=d_{k+i_1}+1=i-1 $, where $d_u'=d_u$ (for $u\ne {k+j_1}, u\ne {k+i_1})$.
Thus,
\begin{eqnarray}
 \nonumber   \frac{\prod_2({T}')}{\prod_2({T}_{max}^2)}  &&=  \frac{(j+1)^{j+1} {(i-1)^{i-1}}}{j^{j} i^{i}}
\\  \nonumber &&=\frac{ \frac{(i-1)^{i-1}} {i^i}}  {\frac {j^{j}} {(j+1)^{(j+1)}}}~\;\;\;\; (\text{by Proposition 2.3})
\\   &&>1.
 \end{eqnarray}

Therefore, we obtain the tree ${T}'$ such that $\prod_2$ is greater than predecessor. So $p=\lfloor \frac{r}{\Delta-2}\rfloor $. In the same way as previously, we also obtain
$\prod_2( T )\leq \Delta^{k\Delta}(\Delta-1)^{p(\Delta-1)} \mu^\mu$.

Finally, this completes our proof.
\end{proof}

\vskip4mm\noindent{\bf Acknowledgements.}


 The work was partially supported by the
National Science Foundation of China under Grant nos.  11271149, 11371162 and 11601006, the Natural
Science Foundation for the Higher Education Institutions of Anhui
Province of China under Grant no. KJ2015A331. Also, it was
partially supported by the Self-determined Research Funds of CCNU
from the colleges basic research and operation of MOE. Furthermore, the authors are grateful to the anonymous referee for a careful checking of the details and for helpful comments that improved this paper.

\end{document}